\documentclass[a4paper,11pt,oneside]{amsart}

\usepackage[utf8]{inputenc}
\usepackage[T1]{fontenc}
\usepackage{amsmath,amssymb,amsfonts,amsthm,mathrsfs}
\usepackage[shortlabels]{enumitem}
\usepackage{comment}
\usepackage{booktabs}
\usepackage{array}
\usepackage{geometry}
\usepackage{multirow}
\usepackage{xcolor}
\usepackage[pagebackref=true,colorlinks=true,linkcolor=blue,citecolor=blue,urlcolor=blue]{hyperref}
\usepackage[alphabetic,backrefs,lite]{amsrefs}
\makeatletter
\let\BR@RerunFileCheck\@empty
\makeatother
\DeclareMathOperator{\codim}{codim}
\DeclareMathOperator{\ord}{ord}
\DeclareMathOperator{\rk}{rk}

\DeclareMathOperator{\Aut}{Aut}

\DeclareMathOperator{\SL}{SL}

\newcommand{\Z}{\mathbb Z}
\newcommand{\F}{\mathbb F}
\newcommand{\PP}{\mathbb P}
\newcommand{\UB}{\mathcal{UB}}
\newcommand{\T}{\mathcal T}

\newcommand{\du}{\mathrm{d}_\mathcal{U}}

\newtheorem{theorem}{Theorem}[section]
\newtheorem{proposition}[theorem]{Proposition}
\newtheorem{lemma}[theorem]{Lemma}

\newtheorem{corollary}[theorem]{Corollary}

\theoremstyle{definition}
\newtheorem{definition}[theorem]{Definition}

\newtheorem{remark}[theorem]{Remark}
\numberwithin{equation}{section}
\hypersetup{pdftitle={On the Beauville Dimension of Finite Abelian Groups}, pdfauthor={Noah Ruhland}}

\begin{document}

\title{On the Beauville Dimension of Finite Abelian Groups}
\author{Noah Ruhland}
\address{University of Bayreuth, Universitätsstr.~30, 95447 Bayreuth, Germany}
\email{noah.ruhland@uni-bayreuth.de}


\begin{abstract}
	Unmixed Beauville manifolds of general type are rigid quotients of products of compact complex curves of higher genus by free diagonal actions of finite groups.
	Their biholomorphism classes are encoded by equivalence classes of group-theoretic data.
	For a finite abelian group $G$, we characterize when these data exist and determine the unmixed Beauville dimension $\du(G)$, the minimum dimension of such a manifold arising from a diagonal action of $G$.
	We prove that $\du(G)<\infty$ if and only if $\exp(G)$ is odd and greater than $3$, and the largest invariant factor of $G$ occurs at least twice.
	In that case,
	\[
		\du(G)=\max\{\rk(G),\rk_3(G)+2\}.
	\]
	Specializing to $\du(G)=2$ recovers the classical result \cite{BCG05}*{Theorem~3.4}: the finite abelian Beauville surface groups are precisely $(\Z/n\Z)^2$ with $n>1$ and $\gcd(n,6)=1$.
	Although the Beauville dimension has a geometric origin, our proof is entirely group-theoretic and combinatorial.
\end{abstract}

\maketitle

\tableofcontents

\section{Introduction}\label{sec:introduction}

An unmixed Beauville surface $X$ is a rigid surface isogenous to a product of curves of unmixed type.
This means that there are compact complex curves $C_1$ and $C_2$ and a finite freely acting group $G \leq \Aut(C_1) \times \Aut(C_2)$ such that
\[
	X \cong (C_1 \times C_2)/G
\]
and that $X$ does not admit any non-trivial deformation.
They were introduced by Catanese in \cite{Cat00} motivated by Beauville's example \cite{B83}*{Exercise~X.13(4)}.

Such surfaces admit a group-theoretic description in terms of two generating triples in the acting group $G$ whose stabilizer sets intersect trivially.
This leads to the following group-theoretic definition.

\begin{definition} \label{def:Beauville-surface-group}
	A finite group $G$ is called an \emph{unmixed Beauville (surface) group} if it admits two triples
	\[
		S_i \in    \mathcal T(G)
		:=
		\bigl\{
		(x,y,z)\in (G\setminus\{1_G\})^3
		\bigm|
		xyz=1,\ \langle x,y,z\rangle=G
		\bigr\}
	\]
	such that
	\[
		\Sigma_{S_1}\cap\Sigma_{S_2}=\{1_G\},
	\]
	where, for $S=(x,y,z)$, the stabilizer set is defined as
	\[
		\Sigma_S
		:=
		\bigcup_{g\in G}
		\left(
		g\langle x\rangle g^{-1}
		\cup g\langle y\rangle g^{-1}
		\cup g\langle z\rangle g^{-1}
		\right).
	\]
	The pair $(S_1,S_2)$ is called an \emph{unmixed Beauville surface structure} on $G$.
\end{definition}

The stabilizer set $\Sigma_{S_i}$ consists precisely of the elements of $G$ fixing a point of the corresponding curve $C_i$.
Hence $\Sigma_{S_1}\cap\Sigma_{S_2}=\{1_G\}$ is equivalent to freeness of the diagonal action on $C_1\times C_2$.
In the surface case, freeness also forces both triples to be hyperbolic, that is, both curves have genus at least two \cite{BCG05}*{Proposition~3.2}.
Hence the resulting surface is of general type.
In the abelian case, these are precisely the groups
\[
	(\mathbb Z/n\mathbb Z)^2 \qquad\text{with}\qquad n>1,\ \gcd(n,6)=1,
\]
see \cite{BCG05}*{Theorem~3.4}. 
The examples go back to Catanese \cite{Cat00}*{Lemma~3.21}.
Determining which non-abelian finite groups are Beauville groups has since motivated extensive research.
Bauer, Catanese and Grunewald proved that the alternating group $\mathcal A_n$ is a Beauville group for all sufficiently large $n$ \cite{BCG05}*{Proposition~3.8}.
Then Garion, Larsen and Lubotzky established the existence of Beauville surface structures for almost all finite simple groups \cite{GLL12}, and Guralnick and Malle completed the simple-group case: every non-abelian finite simple group other than $\mathcal A_5$ is a Beauville group \cite{GM12}.
Fairbairn, Magaard and Parker obtained the corresponding result for finite quasisimple groups, with the exceptions $\mathcal A_5$ and $\SL_2(5)$ \cite{FMP13}.
Beauville structures on finite $p$-groups remain an active research topic \cites{BBF12,FAG17,Fai19,FGGS24}.

The unmixed Beauville dimension studied in this paper arises from extending the theory of Beauville surfaces to higher dimensions.
A first group-theoretic extension was proposed by Carta and Fairbairn in \cite{CF22}.
Motivated by Definition~\ref{def:Beauville-surface-group}, they defined the Beauville dimension of a group $G$ as the smallest integer $r$ for which there exist triples $S_1,\ldots,S_r\in\mathcal T(G)$ satisfying
\[
	\Sigma_{S_1} \cap \dots \cap \Sigma_{S_r} = \{1_G\} .
\]
To account for the geometry in higher dimensions, we also have to allow non-faithful actions on individual curve factors.
This phenomenon already occurs for threefolds and was first described by Frapporti and Gleissner in \cite{FG16}*{Section~3}.
To see how this non-faithfulness enters the group-theoretic data, we will consider the induced actions on the individual curve factors.

For a variety isogenous to a higher product of curves of unmixed type
\[
	X \cong (C_1\times\cdots\times C_r)/G,
\]
the action of $G$ can be restricted to the individual factors $C_i$:
\[
	\psi_i\colon G\longrightarrow \Aut(C_i).
\]
We denote the kernel by $K_i := \ker \psi_i$.
By \cite{FGR26}*{Theorem~2.8}, every such variety admits a minimal realization, unique up to isomorphism and permutation of the factors.
It is characterized by
\[
	\bigcap_{j\ne i}K_j=\{1_G\}
	\qquad\text{for every }i.
\]
For $r=2$, this condition is equivalent to $K_1=K_2=\{1_G\}$, whereas for $r>2$ the individual kernels may be non-trivial.
Consequently, the higher-dimensional Beauville structure associated with a minimal realization records both the kernels $K_i$ and the generating triples $S_i$ in the quotients $G/K_i$ \cite{FGR26}*{Remark~3.3}.

\begin{definition}\label{def:beauville-dimension}
	Let $G$ be a finite group and $r\ge2$.
	An \emph{$r$-fold unmixed Beauville structure} on $G$ is a tuple
	\[
		\bigl[(K_1,S_1),\ldots,(K_r,S_r)\bigr],
	\]
	where $K_i\trianglelefteq G$, $\varphi_i\colon G\longrightarrow G/K_i$ is the quotient map, and $S_i=(x_i,y_i,z_i)\in\mathcal T(G/K_i)$, such that
	\[
		\begin{array}{r@{\quad}l}
			\hypertarget{cond:beauville-hyperbolicity}{\textup{(i)}}
			& S_i \text{ is hyperbolic, i.e.\ } 1 - \frac{1}{\ord(x_i)} - \frac{1}{\ord(y_i)} - \frac{1}{\ord(z_i)} > 0, \quad\text{for every }i,  \\[0.6ex]
			\hypertarget{cond:beauville-minimality}{\textup{(ii)}}
			&
			\displaystyle
			\bigcap_{j\ne i}K_j=\{1_G\}
			\quad\text{for every }i,
			\\[1.4ex]
			\hypertarget{cond:beauville-freeness}{\textup{(iii)}}
			&
			\displaystyle
			\bigcap_{i=1}^r
			\varphi_i^{-1}(\Sigma_{S_i})=\{1_G\}.
		\end{array}
	\]
	We denote the set of $r$-fold unmixed Beauville structures on $G$ by $\UB_r(G)$.
\end{definition}

With the same idea as before, the general definition of the Beauville dimension is as follows.

\begin{definition} \label{def:beauvDim}
	The \emph{(unmixed) Beauville dimension} of a finite group $G$ is defined as
	\[
		\du(G) := \min\{r\ge2\mid\UB_r(G)\ne\emptyset\},
	\]
	with the convention $\du(G)=\infty$ if no such $r$ exists.
\end{definition}

As a first illustration, \cite{FGR26}*{Proposition~5.5} shows that, for $n>1$, $(\Z/n\Z)^3$ admits an unmixed $3$-fold Beauville structure if and only if $\gcd(n,6)=1$.
Thus, in this case, its Beauville dimension is $3$.

Motivated by this result, our goal is to decide which non-trivial finite abelian groups admit unmixed Beauville structures and to compute the Beauville dimension.
To state the result, write
\[
	G\cong \Z/{n_1}\Z \times\cdots\times \Z/{n_s}\Z,
	\qquad 1<n_s\mid\cdots\mid n_1.
\]
We will consider
\[
	\rk_p(G):=\dim_{\mathbb F_p}G[p],
	\qquad
	\rk(G):=\max_p\rk_p(G) = s,
\]
where $G[p]:=\{g\in G\mid pg=0\}$, and let
\[
	\exp(G):=\min \{ n \geq 1 \mid nG = 0 \} = n_1,
	\qquad
	m_{\exp}(G):=
	\bigl|\{i\mid n_i=n_1\}\bigr|
\]
denote the multiplicity of the largest invariant factor.

The main theorem gives both the existence criterion and the exact dimension.

\begin{theorem}\label{thm:main}
	Let $G$ be a finite abelian group.
	Then $\du(G)<\infty$ if and only if $\exp(G)$ is odd and greater than $3$, and $m_{\exp}(G)\ge2$.
	If these conditions hold, then
	\[
		\du(G)=\max\{\rk(G),\rk_3(G)+2\}.
	\]
\end{theorem}

This work is organized as follows.
Section~\ref{sec:definitions} recalls the geometric background and connects it to group theory.
Section~\ref{sec:lower} derives the primary obstructions and lower bounds for the Beauville dimension.
Section~\ref{sec:upper} gives the corresponding constructions to show that the bounds are indeed sharp.
Section~\ref{sec:products} combines the primary components and proves the main theorem, while Section~\ref{sec:examplesAndRemarks} discusses boundary cases and further examples, as well as abelian groups with mixed actions.

\section{From geometry to the Beauville dimension}
\label{sec:definitions}

We first explain how rigidity is encoded by generating triples for the individual factors.
We then assemble the factorwise data and identify the conditions of hyperbolicity, minimality, and freeness.
This leads to the definition of the $r$-fold unmixed Beauville structure and the Beauville dimension.

\subsection{Geometric setup and faithful factor actions}\label{subsec:geometric-setup}

A smooth complex projective variety $X$ is said to be \emph{isogenous to a product of curves of unmixed type} if it admits a presentation
\[
	X\cong (C_1\times\cdots\times C_r)/G,
\]
where the $C_i$ are smooth projective curves of genus at least one and $G\leq\Aut(C_1)\times\cdots\times\Aut(C_r)$ is finite and acts freely and diagonally.
Such a presentation is called an \emph{unmixed realization} of $X$.
The variety $X$ is of general type if and only if $g(C_i)\ge2$ for every $i$.
In this case it is said to be isogenous to a higher product.

The variety $X$ is called \emph{infinitesimally rigid} if
\[
	H^1(X,\Theta_X)=0,
\]
where $\Theta_X$ denotes the holomorphic tangent sheaf.
Throughout this paper, \emph{rigid} always means infinitesimally rigid.
For a free quotient $X=Y/G$, one has
\[
	H^1(X,\Theta_X)\cong H^1(Y,\Theta_Y)^G,
\]
so this condition can be checked on the covering variety by taking $G$-invariants, see \cite{BG21}*{Remark~2.2(2)}.

\begin{definition}
	A rigid variety isogenous to a product of curves is called a \emph{Beauville manifold}.
\end{definition}

The main part of this paper concerns Beauville manifolds of general type with an unmixed realization.
As already explained, for such a realization
\[
	X=(C_1\times\cdots\times C_r)/G,
\]
the diagonal action induces homomorphisms
\[
	\psi_i\colon G\longrightarrow\Aut(C_i).
\]
Set $K_i:=\ker(\psi_i)$ and $Q_i:=G/K_i$.
Thus $Q_i$ acts faithfully on $C_i$.
Since $g(C_i)\ge2$ for every $i$, the Künneth formula gives
\[
	H^1(X,\Theta_X)\cong\bigoplus_{i=1}^r H^1(C_i,\Theta_{C_i})^{Q_i}.
\]
Consequently, $X$ is rigid if and only if $H^1(C_i,\Theta_{C_i})^{Q_i}=0$ for every $i$.
We therefore begin by studying the action on a single factor.

\subsection{Triangle covers and generating triples}\label{subsec:triangle-covers}

Fix an index $i$.
The faithful action of $Q_i$ on $C_i$ defines the quotient map
\[
	f_i\colon C_i\longrightarrow C_i/Q_i.
\]
Using \cite{B83}*{Examples~VI.12(2)}, $H^1(C_i,\Theta_{C_i})^{Q_i}=0$ if and only if $C_i/Q_i\cong\mathbb P^1$ and $f_i$ is branched over exactly three points.
Such a map is called a \emph{triangle cover}, and the corresponding curve $C_i$ is called a \emph{triangle curve}.
Its monodromy is encoded by a generating triple, which we now describe.

Let
\[
	f\colon C\longrightarrow\mathbb P^1
\]
be a connected $H$-Galois triangle cover, where $H$ is finite, and let $B=\{b_1,b_2,b_3\}$ be its branch locus.
Choose a base point $b_0\in\mathbb P^1\setminus B$, a point $q_0\in f^{-1}(b_0)$, and positively oriented loops $\gamma_1,\gamma_2,\gamma_3$ around the branch points such that
\[
	\gamma_1\gamma_2\gamma_3=1.
\]
The resulting monodromy map is a surjection
\[
	\eta_f \colon \pi_1(\mathbb P^1\setminus B,b_0) \longrightarrow H.
\]
Hence
\[
	x:=\eta_f(\gamma_1),\qquad
	y:=\eta_f(\gamma_2),\qquad
	z:=\eta_f(\gamma_3)
\]
satisfy
\[
	xyz=1, \qquad    \langle x,y,z\rangle=H.
\]
Since all three points are branch points, the local monodromies $x,y,z$ are non-trivial.
Thus the marked cover, the chosen point $q_0$, and the loops determine the generating triple $S=(x,y,z)$.
Changing $q_0$ simultaneously conjugates the three entries.

Conversely, if $S=(x,y,z)\in\mathcal T(H)$, then
\[
	\gamma_1\longmapsto x,\qquad
	\gamma_2\longmapsto y,\qquad
	\gamma_3\longmapsto z
\]
defines a surjective monodromy representation.
By the Riemann existence theorem, it determines a connected marked $H$-Galois triangle cover
\[
	C_S\longrightarrow\mathbb P^1
\]
with local monodromies $x,y,z$.
Hence, after fixing the marking and the loops, triangle covers correspond to generating triples, as described in \cite{GT12}*{Section~2} and \cite{FGR26}*{Section~3}.
The triple records the genus of the covering curve and which group elements have fixed points.
The genus is determined by the orders of its entries.
For $S=(x,y,z)$, we define
\[
	\Theta(S)
	=
	1-\frac{1}{\ord(x)}
	-\frac{1}{\ord(y)}
	-\frac{1}{\ord(z)}.
\]
By the Riemann--Hurwitz formula,
\[
	2g(C_S)-2=|H|\Theta(S).
\]
Hence
\[
	g(C_S)\ge2
	\quad\Longleftrightarrow\quad
	\Theta(S)>0.
\]
Such a triple is called \emph{hyperbolic}.

The fixed-point behaviour is recorded by the stabilizer set $\Sigma_S$ introduced in Section~\ref{sec:introduction}; it consists precisely of the elements of $H$ having a fixed point on $C_S$.
If $H$ is abelian, then
\[
	\Sigma_S
	=
	\langle x\rangle
	\cup
	\langle y\rangle
	\cup
	\langle z\rangle.
\]

\subsection{From factor data to a minimal realization}

So far we have encoded one faithful action $Q_i\curvearrowright C_i$.
We now lift these factorwise data to $G$ and assemble them into a Beauville structure of the product.
For fixed $i$, the subgroup
\[
	\bigcap_{j\ne i}K_j
\]
acts trivially on every factor except possibly $C_i$.
The following condition excludes precisely such redundant one-factor actions.

\begin{definition}
	A realization is called \emph{minimal} if
	\[
		\bigcap_{j\ne i}K_j=\{1_G\}
		\qquad\text{for every }i.
	\]
	It is called \emph{absolutely faithful} if $K_i=\{1_G\}$ for every $i$.
\end{definition}

For $r=2$, minimality is equivalent to absolute faithfulness.
For $r\ge3$, minimality is weaker; accordingly, the higher-dimensional formalism must allow generating triples in proper quotients $G/K_i$.

\begin{theorem}[{\cite{FGR26}*{Theorem~2.8}}]
	\label{thm:minimal-realization}
	Every variety isogenous to a higher product of curves of unmixed type admits a minimal realization, unique up to isomorphism and permutation of the factors.
\end{theorem}

In dimension two, Catanese established this result in \cite{Cat00}*{Corollary~3.12 and Proposition~3.13}.

\begin{corollary} \label{cor:realization}
	Let $G$ be a finite group and let $r\ge2$.
	Then $\UB_r(G)\ne\emptyset$ if and only if $G$ occurs as the group of a minimal realization of an $r$-dimensional unmixed Beauville manifold of general type.
	Consequently, $\du(G)$ is the smallest possible dimension of such a manifold.
\end{corollary}

\begin{proof}
	\emergencystretch=1em
	The preceding discussion shows that every minimal realization of an unmixed Beauville manifold of general type gives an element of $\UB_r(G)$.
	Conversely, the Riemann existence theorem constructs the factor curves from an element of $\UB_r(G)$.
	Condition~\hyperlink{cond:beauville-hyperbolicity}{\textup{(i)}} gives general type, conditions~\hyperlink{cond:beauville-minimality}{\textup{(ii)}} and~\hyperlink{cond:beauville-freeness}{\textup{(iii)}} give minimality and freeness, respectively, and rigidity follows from the discussion in Sections~\ref{subsec:geometric-setup} and~\ref{subsec:triangle-covers}.
\end{proof}

\section{Primary obstructions and lower bounds}
\label{sec:lower}

From now on, let $p$ be a prime.
In the abelian cases that follow, groups are written additively.

We first recall the vector space identifications used in the arguments below.

\begin{remark}\label{rem:torsion-reduction}
	Let $G$ be a non-trivial finite abelian $p$-group, and write $n=\rk(G)$, $p^e=\exp(G)$, and $m=m_{\exp}(G)$.
	Then
	\[
		G/pG\cong G[p]\cong\F_p^n,
		\qquad p^{e-1}G\le G[p],
		\qquad p^{e-1}G\cong\F_p^m.
	\]
	All three groups are isomorphic if and only if $m=n$.
	This means that $G$ is homocyclic, i.e.\ $G\cong(\Z/p^e\Z)^n$.
	In this case, multiplication by $p^{e-1}$ induces the canonical isomorphism
	\[
		G/pG\xrightarrow{\sim}G[p]=p^{e-1}G,
		\qquad g+pG\longmapsto p^{e-1}g.
	\]
	For any epimorphism $\varphi\colon G\twoheadrightarrow Q\cong(\Z/p^e\Z)^2$, we have
	\[
		\varphi(G[p])=\varphi(p^{e-1}G)=Q[p].
	\]
	Thus $\varphi|_{G[p]}$ is surjective.
\end{remark}

We now use these identifications to obtain obstructions and lower bounds for the Beauville dimension.
Fix a finite abelian $p$-group $G$ with a $d$-fold unmixed Beauville structure $\bigl[(K_i,S_i)\bigr]_{i=1}^d$.
Write $Q_i:=G/K_i$, with quotient maps $\varphi_i\colon G\to Q_i$.

\subsection{Primary obstructions}

The first obstruction concerns the multiplicity of the largest invariant factor.

\begin{proposition}\label{prop:balanced-necessary}
	Let $G$ be a finite abelian $p$-group.
	\begin{enumerate}[(a)]
		\item If $\du(G)<\infty$, then $m_{\exp}(G)\ge2$.
		\item If $p=2$, then $\du(G)=\infty$.
	\end{enumerate}
\end{proposition}

\begin{proof}
	We first prove (a).
	Since $\du(G)<\infty$, fix a $d$-fold unmixed Beauville structure on $G$ as above.
	Suppose for contradiction that $m_{\exp}(G)=1$, and write $\exp(G)=p^e$.
	Then $L:=p^{e-1}G$ has order $p$.
	Since $\varphi_i$ is surjective,
	\[
		\varphi_i(L)= p^{e-1} \varphi_i(G) =p^{e-1}Q_i.
	\]
	If $\exp(Q_i)<p^e$, then $L\subseteq K_i\subseteq\varphi_i^{-1}(\Sigma_{S_i})$.

	If $\exp(Q_i)=p^e$, then $\varphi_i(L)$ is cyclic of order $p$.
	Since $S_i$ generates $Q_i$, one of its entries $s_i$ has order $p^e$, and hence
	\[
		\varphi_i(L)
		=\langle p^{e-1}s_i\rangle
		\subseteq\Sigma_{S_i}.
	\]
	Thus $L\subseteq\varphi_i^{-1}(\Sigma_{S_i})$ for every $i$.
	Consequently,
	\[
		\{0\} \ne L
		\subseteq\bigcap_{i=1}^d\varphi_i^{-1}(\Sigma_{S_i}),
	\]
	contradicting freeness.

	For part~\textup{(b)}, assume that $p=2$ and write $\exp(G)=2^e$.
	Set $L:=2^{e-1}G\ne\{0\}$.
	Then $\varphi_i(L)=2^{e-1}Q_i$ has $\F_2$-dimension at most two.
	We have already seen that if $\varphi_i(L)$ is cyclic, then $L \subseteq \varphi_i^{-1}(\Sigma_{S_i})$.

	If $\varphi_i(L)$ is not cyclic, then $Q_i\cong(\Z/2^e\Z)^2$.
	Writing $S_i=(x_i,y_i,z_i)$, the classes of $x_i,y_i$ form a basis of $Q_i/2Q_i$.
	Since $\bar z_i=\bar x_i+\bar y_i$, the isomorphism $Q_i/2Q_i\xrightarrow{\sim}2^{e-1}Q_i$ from Remark~\ref{rem:torsion-reduction} maps these three classes to all non-zero elements of $2^{e-1}Q_i$.
	Thus $\varphi_i(L) =2^{e-1}Q_i\subseteq\Sigma_{S_i}$ in every case, contradicting freeness.
\end{proof}

\subsection{Dimension bounds}

We will now use these obstructions to obtain dimension bounds.
Write $n=\rk(G)$ and $p^e=\exp(G)$.
Set
\[
	V:=G[p]\cong\F_p^n,
	\qquad
	\theta_i:=\varphi_i|_V\colon V\longrightarrow Q_i[p].
\]
Then $\dim Q_i[p]=\rk(Q_i)\le2$.

\begin{proposition}\label{prop:general-lower}
	Let $G$ be an abelian $p$-group of rank $n$ that admits a $d$-fold unmixed Beauville structure.
	\begin{enumerate}[(a)]
		\item Then $d \ge n$.
		\item If $p=3$, then $d \ge n + 2$.
	\end{enumerate}
\end{proposition}
\begin{proof}
	For each $i$, choose an entry $s_i$ of $S_i$ of maximal order.
	Since $S_i$ generates $Q_i$, we have $\ord(s_i)=\exp(Q_i)$.
	Set
	\[
		U_i:=\theta_i^{-1}\bigl(\langle s_i\rangle\cap Q_i[p]\bigr).
	\]
	Since $\langle s_i\rangle\cap Q_i[p]$ has codimension at most one in $Q_i[p]$, we have
	\[
		U_i\subseteq\varphi_i^{-1}(\Sigma_{S_i})\cap V
		\qquad\text{and}\qquad
		\codim_V U_i\le1.
	\]
	The freeness condition gives
	\[
		0=\dim\bigcap_{i=1}^dU_i\ge n-d,
	\]
	proving~\textup{(a)}.

	For~\textup{(b)}, assume that $p=3$ and put
	\[
		L:=3^{e-1}G,
		\qquad I:=\{i\mid \dim\theta_i(L)=2\},
		\qquad b:=|I|.
	\]
	For $i\notin I$, our choice of $s_i$ gives
	\[
		\theta_i(L)=3^{e-1}Q_i=\langle3^{e-1}s_i\rangle,
	\]
	so $L\subseteq U_i$.
	Thus
	\[
		H:=\bigcap_{i\notin I}U_i,\qquad h:=\dim H,
	\]
	satisfies $h\ge n-(d-b)$.
	Since $L\subseteq H$, Proposition~\ref{prop:balanced-necessary}(a) also gives $h\ge2$.

	If $i\in I$, then $Q_i\cong(\Z/3^e\Z)^2$.
	Writing $S_i=(x_i,y_i,z_i)$, the elements $3^{e-1}x_i$, $3^{e-1}y_i$, and $3^{e-1}z_i$ determine three of the four lines in $Q_i[3] \cong \F_3^2$.
	Let $M_i=\langle m_i\rangle\cong\F_3$ be the missing line.
	Then
	\[
		Q_i[3]\setminus\Sigma_{S_i}
		=M_i\setminus\{0\}
		=\{m_i,-m_i\}.
	\]
	Since $\dim\theta_i(L)=2$ and $L\subseteq H$, the map $\theta_i|_H$ is surjective.
	Hence
	\[
		H\setminus (\varphi_i^{-1}(\Sigma_{S_i})\cap V) = (\theta_i|_H)^{-1}(\{m_i, -m_i\})
	\]
	is the union of two affine subspaces of dimension $h-2$, neither containing the origin.

	By freeness, the $2b$ affine subspaces obtained in this way cover $H\setminus\{0\}$.
	Jamison's covering theorem \cite{Jam77}*{Theorem~2} says that at least
	\[
		q^{m-k}-1+(q-1)k
	\]
	$k$-dimensional affine subspaces are required to cover $\F_q^m\setminus\{0\}$ without covering the origin.
	Applying this with $q=3$, $m=h$, and $k=h-2$ gives
	\[
		2b\ge 3^2-1+2(h-2)=2h+4,
	\]
	so $b\ge h+2$.
	Since $h\ge n-(d-b)$, we conclude that
	\[
		d=(d-b)+b\ge(n-h)+(h+2)=n+2.
	\]
\end{proof}

The proofs of Propositions~\ref{prop:balanced-necessary} and~\ref{prop:general-lower} use only generation, the product-one relation, and the freeness condition. 
Neither hyperbolicity nor the assumption that the triple entries are non-trivial is needed.

\section{Explicit constructions of unmixed Beauville structures}\label{sec:upper}

We now construct unmixed Beauville structures that attain the lower bounds from Section~\ref{sec:lower}.
The basic idea is to construct linear data on the finite vector space $G[p]$ satisfying minimality and freeness, and then lift them to the whole abelian group $G$.
We treat $p\ge5$ and $p=3$ separately.
The case $p=2$ is excluded by Proposition~\ref{prop:balanced-necessary}(b).

\subsection{Directions and lifting triples}

Fix a prime $p \ge 3$.
The following lemma realizes three arbitrary distinct directions in $\F_p^2$ by the entries of a generating triple.

\begin{lemma}\label{lem:directions}
	Let $L_1,L_2,L_3\le\F_p^2$ be three distinct one-dimensional subspaces.
	Then there exists a triple $(x,y,z)\in\T(\F_p^2)$ such that
	\[
		\langle x\rangle=L_1,\qquad
		\langle y\rangle=L_2,\qquad
		\langle z\rangle=L_3.
	\]
\end{lemma}

\begin{proof}
	Choose non-zero vectors $v_1\in L_1$ and $v_2\in L_2$.
	They form a basis of $\F_p^2$.
	For $v_3\in L_3\setminus\{0\}$, write
	\[
		v_3=av_1+bv_2.
	\]
	Since $L_3$ is distinct from $L_1$ and $L_2$, both $a$ and $b$ are non-zero.
	Set $x:=av_1$, $y:=bv_2$, and $z:=-v_3$.
	Then $x+y+z=0$, the elements $x$ and $y$ form a basis of $\F_p^2$, and $\langle x\rangle=L_1$, $\langle y\rangle=L_2$, and $\langle z\rangle=L_3$.
\end{proof}

The next lemma provides the key reduction step.
It allows us to construct unmixed Beauville structures on $G$ from suitable structures on $G[p]$.

\begin{lemma}\label{lem:lifting}
	Let $G$ be a finite abelian $p$-group of exponent $p^e$ and let $d\ge2$.
	For subgroups $K_1,\ldots,K_d\le G$, define
	\[
		\varphi_i\colon G\longrightarrow Q_i:=G/K_i,
		\qquad
		\theta_i:=\varphi_i|_{G[p]}.
	\]
	Suppose that $Q_i\cong(\Z/p^e\Z)^2$ for every $i$ and that there are triples $\widetilde{S_i}\in\T(Q_i[p])$ such that
	\[
		\bigcap_{j\ne i}\ker\theta_j=\{0\}
		\quad\text{for every }i,
		\qquad
		\bigcap_i\theta_i^{-1}(\Sigma_{\widetilde{S_i}})=\{0\}.
	\]
	Then there exist triples $S_i\in\T(Q_i)$ satisfying $p^{e-1}S_i=\widetilde{S_i}$ such that
	\[
		\bigcap_{j\ne i}K_j=\{0\} \quad\text{for every }i,    \qquad
		\bigcap_{i=1}^d \varphi_i^{-1}(\Sigma_{S_i})=\{0\}.
	\]
	In particular, if $p^e>3$, then
	\[
		\bigl[(K_1,S_1),\ldots,(K_d,S_d)\bigr] \in\UB_d(G).
	\] 
\end{lemma}

\begin{proof}
	Let $\eta_i\colon Q_i/pQ_i\xrightarrow{\sim}Q_i[p]$ be the isomorphism of Remark~\ref{rem:torsion-reduction}.
	Write $\widetilde{S_i}=(\widetilde{x_i},\widetilde{y_i},\widetilde{z_i})$.
	Since $\widetilde{S_i}$ generates $Q_i[p]$ and its entries sum to zero, $\widetilde{x_i},\widetilde{y_i}$ form a basis of $Q_i[p]$.
	Choose representatives $x_i,y_i\in Q_i$ of $\eta_i^{-1}(\widetilde{x_i})$ and $\eta_i^{-1}(\widetilde{y_i})$, respectively.
	Set $z_i:=-x_i-y_i$ and $S_i:=(x_i,y_i,z_i)$.
	Then $p^{e-1}S_i=\widetilde{S_i}$, and the classes of $x_i,y_i$ form a basis of $Q_i/pQ_i$.
	Hence Burnside’s basis theorem \cite{LM02}*{Proposition~1.2.4} yields $Q_i=\langle x_i,y_i\rangle$.
	Moreover, all three entries have order $p^e$, since their images under multiplication by $p^{e-1}$ are non-zero.
	Thus $S_i\in\T(Q_i)$.
	Since every entry of $S_i$ has order $p^e$, we have 
	\[
		\Sigma_{S_i}\cap Q_i[p]=\Sigma_{\widetilde{S_i}}.
	\]
	Any non-zero element of $G$ has a non-zero multiple in $G[p]$, while the kernels and stabilizer sets are closed under taking multiples.
	Therefore, a non-zero element in either $\bigcap_{j\ne i}K_j$ or $\bigcap_i\varphi_i^{-1}(\Sigma_{S_i})$ would yield a non-zero element in the corresponding intersection with $G[p]$.
	The two hypotheses on the maps $\theta_i$ thus imply minimality and freeness, respectively.
	Finally, if $p^e>3$, then
	\[
		\Theta(S_i)=1-\frac{3}{p^e}>0
	\]
	for every $i$.
	Thus the triples are hyperbolic, and
	\[
		\bigl[(K_1,S_1),\ldots,(K_d,S_d)\bigr] \in\UB_d(G).
	\]
\end{proof}

\begin{remark}
	The reduction to finite vector spaces also gives a combinatorial interpretation of minimality and freeness.
	We will not use this interpretation in the proofs below.
	In the setting of Lemma~\ref{lem:lifting}, choose coordinates $V:=G[p]\cong\F_p^n$ and identify $V^*\cong\F_p^n$ using the dual basis.
	For each kernel $E_i:=\ker\theta_i$, its annihilator $E_i^\perp:=\{\lambda\in V^*\mid\lambda|_{E_i}=0\}$ defines a projective line $\ell_i:=\PP(E_i^\perp)$.
	Under the same duality, the three hyperplanes forming $\theta_i^{-1}(\Sigma_{\widetilde{S_i}})$ correspond to three distinct marked points $D_i = \{P_i^{(1)}, P_i^{(2)}, P_i^{(3)}\} \subseteq\ell_i$.
	For the lifted triples, conditions~\textup{(ii)} and~\textup{(iii)} of Definition~\ref{def:beauville-dimension} are respectively equivalent to:
	\[
		\begin{aligned}
			\textup{(ii$'$)}\quad&
			\Bigl\langle\bigcup_{\substack{1\le j\le d\\j\ne i}}\ell_j\Bigr\rangle_{\mathrm{proj}}
			=\PP(\F_p^n)
			&&\forall i\in\{1,\ldots,d\},\\[1ex]
			\textup{(iii$'$)}\quad&
			\langle P_1,\ldots,P_d\rangle_{\mathrm{proj}}=\PP(\F_p^n)
			&&\forall(P_1,\ldots,P_d)\in D_1\times\cdots\times D_d.
		\end{aligned}
	\]
	Here $\langle\cdot\rangle_{\mathrm{proj}}$ denotes the projective span.
\end{remark}

We now construct the quotient maps and the triples for the two cases $p\ge5$ and $p=3$.

\noindent\emph{Common construction of the quotient maps.}
The two constructions below use the same procedure to lift linear maps on $G[p]$ to epimorphisms on $G$.
Let $p\ge3$ be prime, let $e\ge1$ and $r\ge0$, and put $n=r+2$.
For $r\ge1$, choose integers $e\ge f_1\ge\cdots\ge f_r\ge1$, and set
\[
	G=(\Z/p^e\Z)^2\times\prod_{s=1}^r\Z/p^{f_s}\Z,
	\qquad Q=(\Z/p^e\Z)^2.
\]
Write $g=(u,v,w_1,\ldots,w_r)\in G$ and $h=(U,V,W_1,\ldots,W_r)\in\F_p^n$, and use the coordinate isomorphisms
\[
	\begin{aligned}
		\iota_G\colon\F_p^n&\longrightarrow G[p], &
		h&\longmapsto(p^{e-1}U,p^{e-1}V,p^{f_1-1}W_1,\ldots,p^{f_r-1}W_r),\\
		\iota_Q\colon\F_p^2&\longrightarrow Q[p], &
		(X,Y)&\longmapsto(p^{e-1}X,p^{e-1}Y).
	\end{aligned}
\]
Every linear form $\lambda$ on $\F_p^n$ can be written as
\[
	\lambda(h)=\lambda_U U+\lambda_V V+\sum_{s=1}^r\lambda_{W_s}W_s,
	\qquad \lambda_U,\lambda_V,\lambda_{W_s}\in\F_p,
\]
and we call $(\lambda_U,\lambda_V)\in\F_p^2$ its \emph{$(U,V)$-part}.
Choose integer representatives of the coefficients in $\{0,\ldots,p-1\}$ and define
\[
	L_\lambda\colon G\longrightarrow\Z/p^e\Z,
	\qquad L_\lambda(g)=\lambda_U u+\lambda_V v+\sum_{s=1}^r\lambda_{W_s}p^{e-f_s}w_s.
\]
Then $L_\lambda$ is a well-defined homomorphism, and
\[
	L_\lambda(\iota_G(h))=p^{e-1}\lambda(h)
	\qquad\text{for every }h\in\F_p^n.
\]
For each index $i$, choose linear forms $a_i,b_i$ on $\F_p^n$ and define
\[
	\varphi_i\colon G\longrightarrow Q,
	\qquad g\longmapsto\bigl(L_{a_i}(g),L_{b_i}(g)\bigr).
\]
Assume that the $(U,V)$-parts of $a_i$ and $b_i$ are linearly independent.
The matrix of $\varphi_i$ on the first two direct factors then has non-zero determinant modulo $p$.
Hence this restriction is an automorphism of $Q$, and $\varphi_i$ is surjective.
Put $K_i:=\ker\varphi_i$ and identify $Q_i:=G/K_i$ with $Q$ via the induced isomorphism.
Under the coordinate identifications above, the restriction $\theta_i:=\varphi_i|_{G[p]}$ is the linear map
\[
	\theta_i=(a_i,b_i)\colon\F_p^n\longrightarrow\F_p^2.
\]
In particular,
\[
	K_i\cap G[p]=\iota_G\bigl(\ker a_i\cap\ker b_i\bigr).
\]
To apply Lemma~\ref{lem:lifting}, it remains to choose the linear forms $a_i,b_i$ and three distinct directions in each $Q_i[p]$ such that minimality and freeness hold.
By Lemma~\ref{lem:directions}, these directions can be realized by generating triples.

\subsection{The case \texorpdfstring{$p\ge5$}{p >= 5}}

\begin{proposition}\label{prop:pge5-upper}
	Let $p\ge5$ be prime and let $G$ be a finite abelian $p$-group with $m_{\exp}(G)\ge2$.
	Then $\du(G)=\rk(G)$.
\end{proposition}

\begin{proof}
	Let $n:=\rk(G)$ and write $\exp(G)=p^e$.
	If $n=2$, the assertion follows from \cite{Cat00}*{Lemma~3.21}.
	Hence assume $n\ge3$ and put $r=n-2$.

	\smallskip
	\emph{The quotient maps.}
	Put $c:=1$ if $r$ is odd and $c:=2$ if $r$ is even.
	Define $(a_i(h),b_i(h))$ by the following table.
	\renewcommand{\arraystretch}{1.35}
	\[
		\begin{array}{c|c|c|c}
			& n=3,\ p\ge5 & n\ge4,\ p\ge7 & n\ge4,\ p=5\\ \hline
			i=1 & (U-W_1,\,V) & (U,\,V) & (U,\,V)\\
			i=2 & (U,\,V-W_1) & (U,\,V) & (U,\,V)\\
			3\le i\le n-1 & - & (U+W_{i-2},\,V+W_{i-1}) & (U-W_{i-2},\,V-W_{i-1})\\
			i=n & (U-W_1,\,V-W_1) & (U+W_1+W_r,\,V+W_2) & (U-W_r,\,V-cW_1)
		\end{array}
	\]
	Since the $(U,V)$-parts of $a_i$ and $b_i$ are $(1,0)$ and $(0,1)$, the common construction gives epimorphisms $\varphi_i\colon G\to Q$ with $\theta_i=(a_i,b_i)$.

	\smallskip
	\emph{Minimality.}
	We claim that
	\[
		\bigcap_{j\ne i}\ker\theta_j=\{0\}
		\qquad\text{for every }i.
		\tag{M}
	\]
	For $n=3$, any two of the maps in the table have trivial common kernel, so~\textup{(M)} holds.

	Now let $n\ge4$, so $r\ge2$.
	At least one of $\ker\theta_1$ and $\ker\theta_2$ remains in the intersection after omitting the $i$-th kernel.    
	Since $\theta_1(h)=\theta_2(h)=(U,V)$, this forces $U=V=0$.
	For $1\le s\le r-1$, the equations from $\ker\theta_{s+2}$ then reduce to
	\[
		W_s=W_{s+1}=0.
	\]
	The equations from $\ker\theta_n$ reduce to
	\[
		\begin{cases}
			W_1+W_r=0,\quad W_2=0, & p\ge7,\\
			W_r=0,\quad cW_1=0, & p=5.
		\end{cases}
	\]
	The consecutive pairs force all $W_t=0$, except when the first or last pair is omitted.
	In these cases, the equations from $\ker\theta_n$ force the remaining coordinates to vanish, proving~\textup{(M)}.

	\smallskip
	\emph{The slope sets.}
	We use the following convention.
	A finite slope $\mu\in\F_p$ denotes the line $\langle(1,\mu)\rangle\le\F_p^2$, while $\infty$ denotes the vertical line $\langle(0,1)\rangle\le\F_p^2$.
	Choose triples of lines in $Q[p]\cong\F_p^2$ whose slope sets are given as follows.
	\[
		\begin{array}{c|c|w{c}{2cm}|w{c}{2cm}}
			&  & n=3 & n\ge4 \\[0.5ex] \hline
			\multirow{4}{*}{$p>7$}
			& \rule{0pt}{2.4ex}i=1 & \{-1,0,1\} & \{-1,0,1\}\\
			& i=2 & \{-2,0,2\} & \{2,3,4\}\\
			& 3\le i\le n-1 & - & \{-1,0,1\}\\
			& i=n & \{-1,1,3\} & \{2,3,4\}\\[0.5ex]  \hline
			\multirow{5}{*}{$p=7$}
			& \rule{0pt}{2.4ex}i=1 & \{-1,2,3\} & \{-1,0,1\}\\
			& i=2 & \{0,1,2\} & \{2,3,4\}\\
			& i=3 & - & \{1,4,5\}\\
			& 4\le i\le n-1 & - & \{-1,0,1\}\\
			& i=n & \{-1,1,4\} & \{-1,0,3\}\\[0.5ex] \hline
			\multirow{4}{*}{$p=5$}
			& \rule{0pt}{2.4ex}i=1 & \{-1,1,2\} & \{-1,0,1\}\\
			& i=2 & \{0,1,3\} & \{2,3,\infty\}\\
			& 3\le i\le n-1 & - & \{2,3,\infty\}\\
			& i=n & \{2,3,\infty\} & \{2,3,\infty\}
		\end{array}
	\]
	Each slope set in the table consists of three distinct points of $\PP^1(\F_p)$.
	By Lemma~\ref{lem:directions}, choose a generating triple $\bar S_i\in\T(Q_i[p])$ with precisely these three directions.

	\smallskip
	\emph{Freeness.}
	Suppose that $h\in\F_p^n\cong G[p]$ lies in every stabilizer preimage.
	For each $i$, choose a slope $\mu_i$ from the corresponding table entry such that $\theta_i(h)$ lies on that line.
	The corresponding equation is
	\[
		\begin{cases}
			(b_i-\mu_i a_i)(h)=0,&\text{if } \mu_i\ne\infty,\\
			a_i(h)=0,&\text{if } \mu_i=\infty.
		\end{cases}
	\]
	Let $M(\mu_1,\ldots,\mu_n)$ be the coefficient matrix of these equations, with columns ordered as $(V,U,W_1,\ldots,W_r)$.
	It suffices to prove that $M(\mu_1,\ldots,\mu_n)$ is invertible for every permitted choice of slopes.

	For $n=3$ and finite $\mu_3$, its determinant is
	\[
		E_3(\mu_1,\mu_2,\mu_3)
		= \det \begin{pmatrix}
			1 & -\mu_1 & \mu_1\\
			1 & -\mu_2 & -1\\
			1 & -\mu_3 & \mu_3-1
		\end{pmatrix}
		=(\mu_1+1)(\mu_2-\mu_3)-\mu_3(\mu_2-\mu_1).
	\]
	The only vertical case occurs for $p=5$ and gives
	\[
		E_3(\mu_1,\mu_2,\infty)
		=\det \begin{pmatrix}
			1 & -\mu_1 & \mu_1\\
			1 & -\mu_2 & -1\\
			0 & 1 & -1
		\end{pmatrix}
		=\mu_2+1.
	\]
	For completeness, direct substitution in the finite formula gives
	\[
		\begin{array}{c|c}
			 &\text{possible values of }E_3\\ \hline
			p>7&\{-9,-7,-5,-3,-1,1,3,5,7\}\subseteq\F_p^\times,\\
			p=7&\{1,2,3,4,5,6\}\subseteq\F_7^\times,\\
			p=5&\{1,2,3,4\}\subseteq\F_5^\times,
		\end{array}
	\]
	and in the vertical case
	\[
		\mu_2+1\in\{1,2,4\}\subseteq\F_5^\times.
	\]
	Hence $E_3\ne0$ in all cases.

	For $n\ge4$ and $p\ne5$, the chosen linear forms give
	\begin{align*}
		E_n(\mu_1,\ldots,\mu_n)
		&= \det \begin{pmatrix}
			1 & -\mu_1 & 0 & 0 & 0 & \cdots & 0 & 0\\
			1 & -\mu_2 & 0 & 0 & 0 & \cdots & 0 & 0\\
			1 & -\mu_3 & -\mu_3 & 1 & 0 & \cdots & 0 & 0\\
			1 & -\mu_4 & 0 & -\mu_4 & 1 & \cdots & 0 & 0\\
			1 & -\mu_5 & 0 & 0 & -\mu_5 & \ddots & 0 & 0\\
			\vdots & \vdots & \vdots & \vdots & \ddots & \ddots & 1 & 0\\
			1 & -\mu_{n-1} & 0 & 0 & 0 & \cdots & -\mu_{n-1} & 1\\
			1 & -\mu_n & -\mu_n & 1 & 0 & \cdots & 0 & -\mu_n
		\end{pmatrix} = \\
		&= (-1)^{n-1}(\mu_2-\mu_1)
		\left( \mu_n-\mu_3+\prod_{k=3}^{n}\mu_k \right).
		\tag{D}
	\end{align*}
	For $p>7$, the factor $\mu_2-\mu_1$ is non-zero.
	If one of $\mu_3,\ldots,\mu_{n-1}$ is zero, the second factor in~\textup{(D)} is $\mu_n-\mu_3\ne0$.
	Otherwise put
	\[
		\varepsilon:=\prod_{k=3}^{n-1}\mu_k\in\{1,-1\}.
	\]
	If $\varepsilon=-1$, the second factor equals $-\mu_3\in\{1,-1\}$; if $\varepsilon=1$, it equals
	\[
		2\mu_n-\mu_3\in\{3,5,7,9\}.
	\]
	These values are non-zero modulo every prime $p>7$.

	For $p=7$, the sets for $\mu_1$ and $\mu_2$ are disjoint, so the first factor in~\textup{(D)} is again non-zero.
	If $\mu_n=0$ or some $\mu_k=0$ with $4\le k\le n-1$, the second factor is $\mu_n-\mu_3\ne0$, since
	\[
		\mu_3\in\{1,4,5\},
		\qquad
		\mu_n\in\{-1,0,3\}.
	\]
	Otherwise set $ \varepsilon=\prod_{k=4}^{n-1}\mu_k\in\{1,-1\}$.
	The second factor becomes
	\[
		\mu_n-\mu_3+\varepsilon\mu_3\mu_n.
	\]
	For $\varepsilon=-1$ its possible values are $\{-1,1,4\}$, and for $\varepsilon=1$ they are $\{-1,3,4,5\}$.
	Thus $E_n\ne0$ also for $p=7$.

	It remains to treat $p=5$ and $n\ge4$.
	Recall that
	\[
		\mu_1\in\{-1,0,1\},\qquad
		\mu_i\in\{2,3,\infty\}\quad \text{for } \ 2\le i\le n.
	\]
	The first two rows of $M$ have zero $(W_1,\dots,W_r)$-part.
	Thus, we can write
	\[
		M = \begin{pmatrix}
			A_0 & 0\\
			* & B
		\end{pmatrix}.
	\]
	Its upper-left block is
	\[
		A_0=
		\begin{cases}
			\begin{pmatrix}1&-\mu_1\\1&-\mu_2\end{pmatrix},
			&\mu_2\ne\infty,\\[2ex]
			\begin{pmatrix}1&-\mu_1\\0&1\end{pmatrix},
			&\mu_2=\infty.
		\end{cases}
	\]
	Consequently,
	\[
		\det A_0=
		\begin{cases}
			\mu_1-\mu_2,&\mu_2\ne\infty,\\
			1,&\mu_2=\infty,
		\end{cases}
	\]
	which is non-zero because the finite slope sets $\{-1,0,1\}$ and $\{2,3\}$ are disjoint in $\F_5$.
	We now focus on the lower-right $r\times r$ block $B$.
	For $3\le s\le n$ put
	\[
		(\varepsilon_s,\eta_s)=
		\begin{cases}
			(\mu_s,1),&\mu_s\ne\infty,\\
			(-1,0),&\mu_s=\infty.
		\end{cases}
	\]
	Then
	\[
		B=
		\begin{pmatrix}
			\varepsilon_3&-\eta_3&&&0\\
			0&\varepsilon_4&-\eta_4&&\\
			\vdots&&\ddots&\ddots&\\
			0&&&\varepsilon_{n-1}&-\eta_{n-1}\\
			-c\eta_n&0&\cdots&0&\varepsilon_n
		\end{pmatrix}.
	\]
	Expanding the determinant along the first column gives
	\[
		\det B = \varepsilon_3\prod_{s=4}^n\varepsilon_s +
		(-1)^{r+1} (-c\eta_n) \prod_{s=3}^{n-1}(-\eta_s) =
		\prod_{s=3}^n\varepsilon_s -c\prod_{s=3}^n\eta_s.
	\]
	If at least one $\mu_s$ is vertical, then $\prod_s\eta_s=0$, whereas every $\varepsilon_s$ is non-zero, and hence $\det B\ne0$.
	If all slopes $\mu_s$ are finite, then
	\[
		\det B=\prod_{s=3}^n\mu_s-c.
	\]
	The product lies in $\{2,3\}$ when $r$ is odd and in $\{1,4\}$ when $r$ is even.
	By definition, it therefore differs from $c$.
	Hence $\det B\ne0$, and
	\[
		\det M=(\det A_0)(\det B)\ne0.
	\]
	In summary, we have proved the freeness condition
	\[
		\bigcap_{i=1}^n
		\theta_i^{-1}(\Sigma_{\bar S_i})=\{0\}
		\tag{F}
	\]
	for every prime $p\ge5$ and $n\ge3$.
	Together with~\textup{(M)}, the identifications $G[p]\cong \F_p^n$, $Q_i\cong Q$, and $p^e>3$, Lemma~\ref{lem:lifting} yields an $n$-fold unmixed Beauville structure on $G$.
	Proposition~\ref{prop:general-lower}(a) supplies the opposite inequality, so $\du(G)=n$.
\end{proof}

\subsection{The case \texorpdfstring{$p=3$}{p = 3}}

Hyperbolicity requires exponent at least $3^2$, since triples in groups of exponent $3$ have $\Theta(S)\le0$.
We also record the elliptic construction at exponent $3$, which is used below.

\begin{proposition}\label{prop:p3-upper}
	Let $G$ be a non-trivial finite abelian $3$-group with $m_{\exp}(G)\ge2$, and put $n:=\rk(G)$.
	If $\exp(G)\ge9$, then
	\[
		\du(G)=n+2.
	\]
\end{proposition}

\begin{proof}
	Write $\exp(G)=3^e$ and choose a decomposition of $G$ as in the common construction, with $r=n-2$.
	Put $m=n+2$ and define linear forms on $\F_3^n\cong G[3]$ by
	\[
		\lambda_1(h)=\lambda_{m+1}(h):=U,
		\qquad \lambda_2(h):=V,
		\qquad \lambda_3(h):=U+V,
		\qquad \lambda_4(h):=U-V,
	\]
	and, for $1\le t\le r$,
	\[
		\lambda_{4+t}(h):=\begin{cases}
			V+W_t,&t\text{ odd},\\
			U+V+W_t,&t\text{ even}.
		\end{cases}
	\]

	\smallskip
	\emph{The quotient maps.}
	For $1\le i\le m$, the $(U,V)$-parts of $\lambda_i$ and $\lambda_{i+1}$ are linearly independent.
	Hence the common construction with $a_i=\lambda_i$ and $b_i=\lambda_{i+1}$ gives epimorphisms $\varphi_i\colon G\to Q$ with
	\[
		\theta_i=(\lambda_i,\lambda_{i+1})\colon\F_3^n\longrightarrow\F_3^2.
	\]

	\smallskip
	\emph{Minimality.}
	Fix $i$.
	Since the pairs $(\lambda_j,\lambda_{j+1})$ occur cyclically, an element of $\bigcap_{j\ne i}\ker\theta_j$ lies in the kernel of every $\lambda_k$.
	The equations $\lambda_1(h)=\lambda_2(h)=0$ give $U=V=0$, and then $\lambda_{4+t}(h)=0$ gives $W_t=0$ for every $t$.
	Hence
	\[
		\bigcap_{j\ne i}\ker\theta_j=\{0\}. \tag{M$_3$}\label{eq:p3-minimality}
	\]

	\smallskip
	\emph{Freeness.}
	By Lemma~\ref{lem:directions}, choose in each copy of $\F_3^2$ a triple $\bar S_i$ with directions
	\[
		\langle(1,0)\rangle,
		\langle(1,1)\rangle,
		\langle(1,-1)\rangle.
	\]
	Its stabilizer set is the union of these three lines and therefore misses the non-zero points on the vertical line $\langle(0,1)\rangle$.

	Let $0\ne h\in G[3]$.
	At least one of
	\[
		\lambda_1(h),\lambda_2(h),\lambda_3(h),\lambda_4(h)
	\]
	is zero, since for every $(U,V)\in\F_3^2$ one of
	\[
		U,\qquad V,\qquad U+V,\qquad U-V
	\]
	vanishes.
	On the other hand, not all $\lambda_k(h)$ vanish: otherwise $\lambda_1(h)=\lambda_2(h)=0$ would give $U=V=0$, and the forms $\lambda_{4+t}(h)$ would then give $W_t=0$ for every $t$, contradicting $h\ne0$.
	Hence the cyclic sequence $\lambda_1(h),\ldots,\lambda_m(h)$ contains a zero followed by a non-zero entry, say
	\[
		\lambda_i(h)=0,
		\qquad \lambda_{i+1}(h)\ne0.
	\]
	Then $\theta_i(h)$ is a non-zero vertical vector and is not contained in $\Sigma_{\bar S_i}$.
	Hence
	\[
		\bigcap_{i=1}^m\theta_i^{-1}(\Sigma_{\bar S_i})=\{0\}. 
		\tag{F$_3$} 
		\label{eq:p3-freeness}
	\]
	If $e\ge2$, conditions~\textup{(M$_3$)} and~\textup{(F$_3$)} allow us to apply Lemma~\ref{lem:lifting}.
	In this case, we obtain an $(n+2)$-fold unmixed Beauville structure on $G$.
	Proposition~\ref{prop:general-lower}(b) gives the reverse inequality, and hence $\du(G)=n+2$.
\end{proof}

The construction also applies when $e=1$ and gives $n+2$ pairs satisfying minimality and freeness, with generating triples of type $(3,3,3)$.

\section{From primary components to the main theorem}
\label{sec:products}

We now assemble the primary constructions.
First, we combine structures on groups of coprime order.
Then we apply this construction to the Sylow decomposition of a finite abelian group.

\subsection{Coprime products}

We begin by recording how Beauville structures combine across factors of coprime order.

\begin{proposition}\label{prop:coprime-product}
	Let $N_1,N_2$ be finite groups with $\gcd(|N_1|,|N_2|)=1$.
	Then
	\[
		\du(N_1\times N_2)
		\le \max\{\du(N_1),\du(N_2)\}.
	\]
\end{proposition}

\begin{proof}
	If one of the two dimensions is infinite, there is nothing to prove.
	Put
	\[
		d:=\max\{\du(N_1),\du(N_2)\}.
	\]
	Choose $d$-fold unmixed Beauville structures on $N_1$ and $N_2$, repeating pairs if necessary, and write them as
	\[
		\bigl[(K_i^{(\nu)},S_i^{(\nu)})\bigr]_{i=1}^d,
		\qquad \nu=1,2,
	\]
	where
	\[
		S_i^{(\nu)}
		=\bigl(x_i^{(\nu)},y_i^{(\nu)},z_i^{(\nu)}\bigr) \in \T(N_\nu/K_i^{(\nu)}).
	\]
	Set
	\[
		K_i:=K_i^{(1)}\times K_i^{(2)}
	\]
	and
	\[
		S_i:=\bigl(
		(x_i^{(1)},x_i^{(2)}),
		(y_i^{(1)},y_i^{(2)}),
		(z_i^{(1)},z_i^{(2)})
		\bigr).
	\]
	Denote the corresponding quotient maps by $\varphi_i^{(\nu)}$, and, under the natural identification of the quotients, write $\varphi_i=\varphi_i^{(1)}\times\varphi_i^{(2)}$.

	The product-one relation holds componentwise.
	Since $N_1$ and $N_2$ have coprime orders,
	\[
		\langle(u_1,u_2)\rangle
		=\langle u_1\rangle\times\langle u_2\rangle
	\]
	for corresponding entries $u_1,u_2$.
	Thus $S_i$ generates
	\[
		(N_1\times N_2)/K_i
		\cong (N_1/K_i^{(1)})\times(N_2/K_i^{(2)}).
	\]
	The orders of corresponding entries multiply, and hence
	\[
		\Theta(S_i)\ge\Theta(S_i^{(1)})>0.
	\]

	The minimality condition holds componentwise:
	\[
		\bigcap_{j\ne i}K_j
		=
		\left(\bigcap_{j\ne i}K_j^{(1)}\right)
		\times
		\left(\bigcap_{j\ne i}K_j^{(2)}\right)
		=\{1\}.
	\]
	Finally,
	\[
		\Sigma_{S_i}
		\subseteq
		\Sigma_{S_i^{(1)}}\times\Sigma_{S_i^{(2)}}.
	\]
	Therefore
	\[
		\bigcap_{i=1}^d\varphi_i^{-1}(\Sigma_{S_i})
		\subseteq
		\left(
		\bigcap_{i=1}^d
		(\varphi_i^{(1)})^{-1}(\Sigma_{S_i^{(1)}})
		\right) \,
		\times \,
		\left(
		\bigcap_{i=1}^d
		(\varphi_i^{(2)})^{-1}(\Sigma_{S_i^{(2)}})
		\right)
		=\{1\}.
	\]
	Thus the pairs $(K_i,S_i)$ form a $d$-fold unmixed Beauville structure on $N_1\times N_2$.
\end{proof}

\begin{remark}\label{rem:coprime-product-equality}
	\begin{enumerate}[(a)]
		\item \emph{Projection to coprime factors.}
			Let $[(K_i,S_i)]_{i=1}^d\in\UB_d(N_1\times N_2)$, with $\gcd(|N_1|,|N_2|)=1$.
			We now reverse the preceding construction by projecting the data to the two factors.
			Coprimality gives
			\[
				K_i=K_i^{(1)}\times K_i^{(2)},
				\qquad K_i^{(\nu)}:=K_i\cap N_\nu.
			\]
			Put $Q_i^{(\nu)}:=N_\nu/K_i^{(\nu)}$ and let $S_i^{(\nu)}$ be the projection of $S_i$ to $Q_i^{(\nu)}$.
			These data still satisfy the product-one, generation, minimality, and freeness conditions.
			Projection may nevertheless produce trivial entries or destroy hyperbolicity.
			If $d_0:=\du(N_1\times N_2)<\infty$, equality holds in Proposition~\ref{prop:coprime-product} if and only if there is an unmixed Beauville structure in $\UB_{d_0}(N_1\times N_2)$ that projects to hyperbolic triples on both factors.

		\item \emph{Removing trivial projected entries.}
			Let $[(K_i,S_i)]_{i=1}^d\in\UB_d(G)$, where $G$ is finite abelian of odd order.
			We modify this structure so that every Sylow projection has non-zero entries, showing that at least this problem from part~(a) is avoidable.
			The idea is to keep $d$ fixed and shrink kernels if necessary.
			For each $p\mid |G|$, write
			\[
				K_{i,p}:=K_i\cap G_p,  \qquad Q_{i,p}:=G_p/K_{i,p}
			\]
			and let $S_{i,p}$ be the projected triple. \\ 
			First, correct the triples in non-trivial quotients.
			A zero entry forces $S_{i,p}=(0,a,-a)$ up to permutation and $Q_{i,p}=\langle a\rangle$.
			Keeping $K_{i,p}$ fixed, replace this triple by $(a,a,-2a)$ in the corresponding order.
			Since $p$ is odd, the stabilizer set remains $Q_{i,p}$ and no entry order decreases. \\ 
			Next, replace the trivial quotients.
			Choose $j$ with $Q_{j,p}\ne0$, which exists by freeness.
			Whenever $Q_{i,p}=0$, replace the pair
			\[
				(K_{i,p},S_{i,p})=(G_p,(0,0,0))
				\quad\text{by}\quad (K_{j,p},S_{j,p}),
			\]
			using a triple already corrected in the first step.
			This reduces the kernel and does not enlarge the stabilizer preimage, which was $G_p$. \\
			Finally, recombining the modified data as in the proof of Proposition~\ref{prop:coprime-product} preserves minimality, freeness, and hyperbolicity, since kernels and stabilizer preimages have not increased and entry orders have not decreased.
	\end{enumerate}
\end{remark}

For a generating triple $S=(x,y,z)$, its \emph{type} is $(\ord(x),\ord(y),\ord(z))$, considered up to permutation.
The following remark shows that the projected-hyperbolicity hypothesis in Remark~\ref{rem:coprime-product-equality}(a) is essential.

\begin{remark}\label{rem:hyperbolicity-lost}
	Consider
	\[
		G=(\Z/15\Z)^2\cong\F_3^2\times\F_5^2.
	\]
	Each generating triple in $\F_3^2$ has type $(3,3,3)$.
	We choose four such triples as in the proof of Proposition~\ref{prop:p3-upper}, so that the freeness condition holds.
	We combine them with the two triples of a Beauville structure on $\F_5^2$, using each twice.
	With trivial kernels, these four triples $S_i$ satisfy minimality and freeness and have type $(15,15,15)$.
	In particular,
	\[
		\Theta(S_i)=1-\frac{3}{15}>0.
	\]
	By Theorem~\ref{thm:main}, we have $\du(G)=4$, but the projected triples on $\F_3^2$ have $\Theta=0$ and $\du(\F_3^2)=\infty$.
	Thus projection to a Sylow subgroup does not always preserve hyperbolicity.
\end{remark}

\subsection{Proof of the main theorem}

We now apply the preceding results to the Sylow decomposition of a finite abelian group.

\begin{proof}[Proof of Theorem~\ref{thm:main}]
	Write
	\[
		G=\prod_{p\mid |G|}G_p,
		\qquad
		n_p:=\rk_p(G),
	\]
	where $G_p$ denotes the Sylow $p$-subgroup of $G$.
	Suppose first that $\du(G)<\infty$, put $d:=\du(G)$, and choose
	\[
		\bigl[(K_1,S_1),\ldots,(K_d,S_d)\bigr]
		\in\UB_d(G).
	\]
	For each $i$ and $p$, put $K_{i,p}:=K_i\cap G_p$.
	Under the decomposition
	\[
		G/K_i\cong\prod_{p\mid |G|}G_p/K_{i,p},
	\]
	let $S_{i,p}$ be the projection of $S_i$ to $G_p/K_{i,p}$.
	By Remark~\ref{rem:coprime-product-equality}(a), for each fixed $p$, the pairs $(K_{i,p},S_{i,p})$ still satisfy the generation, minimality, and freeness conditions, although the triples need not be hyperbolic.
	The obstruction in Proposition~\ref{prop:balanced-necessary}(b) yields $G_2=0$.
	By Remark~\ref{rem:coprime-product-equality}(b), we may then assume that all projected entries are non-trivial, without changing $d$.
	By the note after Proposition~\ref{prop:general-lower}, we obtain
	\[
		m_{\exp}(G_p)\ge2 \qquad\text{for every }p\mid |G|,
	\]
	which is equivalent to $m_{\exp}(G)\ge2$.
	Moreover, $\exp(G)\ne3$, since a generating triple in a group of exponent $3$ cannot be hyperbolic.
	Hence $\exp(G)>3$ is odd.

	By Proposition~\ref{prop:general-lower}(a), we obtain
	\[
		\du(G)\ge\max_p n_p=\rk(G),
	\]
	and part~\textup{(b)} yields, if $3\mid |G|$,
	\[
		\du(G)\ge n_3+2.
	\]
	If $3\nmid |G|$, then $n_3=0$ and $\rk(G)\ge2$.
	Therefore, in all cases,
	\[
		\du(G)
		\ge\max\{\rk(G),\rk_3(G)+2\}.
		\tag{$\ast$}
	\]
	Conversely, suppose that $\exp(G)>3$ is odd and $m_{\exp}(G)\ge2$, and put
	\[
		D:=\max\{\rk(G),\rk_3(G)+2\}.
	\]
	Then $m_{\exp}(G_p)\ge2$ for every Sylow subgroup $G_p$.
	Proposition~\ref{prop:pge5-upper} provides an $n_p$-fold structure on $G_p$ for every prime $p\ge5$ dividing $|G|$.
	If $3\mid |G|$, Proposition~\ref{prop:p3-upper} provides $n_3+2$ pairs satisfying generation, minimality, and freeness.
	The triples are hyperbolic if and only if $\exp(G_3)>3$.
	Repeating pairs if necessary, we combine these families componentwise as in Proposition~\ref{prop:coprime-product}, preserving generation, minimality, and freeness.
	Since $\exp(G)>3$, at least one primary family is hyperbolic, and hence so are all product triples.
	Thus $G$ admits a $D$-fold unmixed Beauville structure, so
	\[
		\du(G)\le D.
	\]
	Together with $(\ast)$, this proves
	\[
		\du(G)=\max\{\rk(G),\rk_3(G)+2\}.
	\]
\end{proof}

\section{Boundary cases and further examples} \label{sec:examplesAndRemarks}

This final section clarifies the scope of the preceding results.
We record families of infinite Beauville dimension.
We then consider what happens when we drop the hyperbolicity condition or allow mixed actions.

\subsection{Groups of infinite Beauville dimension}

As noted in Section~\ref{sec:introduction}, the only finite quasisimple groups that are not Beauville surface groups are the alternating group $\mathcal A_5$ and the special linear group $\SL_2(5)$.
We therefore determine their Beauville dimensions, as well as those of the dihedral groups.

The proofs of Propositions~\ref{prop:A5andSL25} and~\ref{prop:dihedral} adapt the surface-case arguments in \cite{FJ11}*{Theorem~5.1, Case~8} and \cite{BCG05}*{Lemma~3.7}, respectively.

\begin{proposition} \label{prop:A5andSL25}
	The groups $\mathcal A_5$ and $\SL_2(5)$ have infinite Beauville dimension.
\end{proposition}
\begin{proof}
	First, suppose that
	\[
		\bigl[(K_1,S_1),\ldots,(K_d,S_d)\bigr]
		\in\mathcal{UB}_d(\mathcal A_5).
	\]
	Since $\mathcal A_5$ is simple, every $K_i$ is trivial.
	Write $S_i=(x_i,y_i,z_i)$.
	We first observe that every generating triple of $\mathcal A_5$ contains an element of order $5$.
	Otherwise, the orders of all three entries lie in $\{2,3\}$.
	A generating triple of type $(a,b,c)$ induces an epimorphism
	\[
		\Delta(a,b,c)
		:=
		\langle
		\alpha,\beta,\gamma
		\mid
		\alpha^a=\beta^b=\gamma^c=\alpha\beta\gamma=1
		\rangle
		\longrightarrow\mathcal A_5.
	\]
	Up to permutation, the only possibilities are
	\[
		\begin{alignedat}{2}
			\Delta(2,2,2) &\cong (\mathbb Z/2\mathbb Z)^2, \qquad\qquad&
			\Delta(2,2,3) &\cong \mathfrak S_3, \\
			\Delta(2,3,3) &\cong \mathcal A_4, &
			\Delta(3,3,3) &\cong \mathbb Z^2\rtimes \mathbb Z/3\mathbb Z .
		\end{alignedat}
	\]
	All these groups are solvable.
	Hence every quotient of any of them is solvable, whereas $\mathcal A_5$ is not, a contradiction.

	Since all subgroups of order $5$ in $\mathcal A_5$ are conjugate, it follows that
	\[
		\{g\in\mathcal A_5\mid \ord(g)=5\}
		\subseteq
		\bigcap_{i=1}^d \Sigma_{S_i},
	\]
	contradicting freeness.

	Second, let $G=\SL_2(5)$ and suppose that
	\[
		\bigl[(K_1,S_1),\ldots,(K_d,S_d)\bigr]
		\in\mathcal{UB}_d(G).
	\]
	Since $G$ is quasisimple, every proper normal subgroup is contained in $Z(G)$.
	Thus $K_i\le Z(G)$ for all $i$, and we can consider the quotient map
	\[
		G/K_i\longrightarrow G/Z(G)\cong\mathcal A_5.
	\]
	The image of $S_i$ is a generating triple of $\mathcal A_5$.
    By the preceding argument, one of its entries has order $5$.
	Hence, $\varphi_i^{-1}(\Sigma_{S_i})$ contains an order-$5$ subgroup of $G$.
	Since all subgroups of order $5$ in $G$ are conjugate, it follows that
	\[
		\{g\in G\mid \ord(g)=5\}
		\subseteq
		\bigcap_{i=1}^d \varphi_i^{-1}(\Sigma_{S_i}),
	\]
	again contradicting freeness.
\end{proof}

\begin{proposition} \label{prop:dihedral}
	Every dihedral group $D_{2m}$ has infinite Beauville dimension.
\end{proposition}
\begin{proof}
	Let
	\[
		G=D_{2m}=\langle r,s\mid r^m=s^2=1,\ srs=r^{-1}\rangle,
		\qquad m\ge2 .
	\]
	Suppose that
	\[
		\bigl[(K_1,S_1),\ldots,(K_d,S_d)\bigr]
		\in\mathcal{UB}_d(G),
	\]
	and put $R=\langle r\rangle$.
	If $K_i\not\le R$, then $K_i$ contains a reflection.
	Hence the image of $r$ in $G/K_i$ has order at most $2$, and therefore $G/K_i$ is abelian and every element has order at most $2$, contradicting the hyperbolicity condition $\Theta(S_i)>0$.
	Thus $K_i\le R$ for every $i$.

	Set $R_i=R/K_i$.
	Then $G/K_i$ is dihedral with rotation subgroup $R_i$.
	For $S_i=(x_i,y_i,z_i)$, exactly one entry lies in $R_i$: not all three can lie in $R_i$, and the product relation forces an even number of entries outside $R_i$.
	Let this entry be $x_i$.
	Since $S_i$ generates $G/K_i$ with $y_i z_i = x_i^{-1}$, we have $\langle x_i\rangle=R_i$.
	Consequently,
	\[
		R_i\subseteq\Sigma_{S_i}, \qquad \text{and hence} \qquad R\subseteq \varphi_i^{-1}(\Sigma_{S_i}) \ \text{ for every } i.
	\]
	Therefore
	\[
		\{1\} \neq R\subseteq\bigcap_{i=1}^d \varphi_i^{-1}(\Sigma_{S_i}),
	\]
	contradicting freeness.

	Alternatively, we could have argued that a quotient with $K_i\not\le R$ has exponent at most two.
	If $K_i\le R$, every generating product-one triple has type $(2,2,|R/K_i|)$ and hence $\Theta(S_i)=-1/|R/K_i|<0$.
	Thus no quotient admits a hyperbolic generating triple.
\end{proof}

\subsection{The elliptic 3-primary case} \label{subsec:specialTypeAndMixed}

We now drop the hyperbolicity condition~\textup{(i)} of Definition~\ref{def:beauville-dimension} and consider $G=(\Z/3\Z)^n$ with $n\ge2$.
The triples of type $(3,3,3)$ from Proposition~\ref{prop:p3-upper} correspond to elliptic triangle covers, which we use to construct smooth infinitesimally rigid quotients.

\begin{proposition}\label{prop:EllipticRigid}
	Let $E=\mathbb C/\mathbb Z[\zeta_3]$, where $\zeta_3=\exp(\frac{2\pi i}{3})$.
	For every $n\ge2$, the group $G=(\Z/3\Z)^n$ admits a minimal free diagonal action on $E^{n+2}$ whose quotient is infinitesimally rigid and has Kodaira dimension~$0$.
\end{proposition}

\begin{proof}
	The point $t:=(1-\zeta_3)^{-1} \in E$ has order $3$ and satisfies $\zeta_3t=t$, yielding a faithful action
	\[
		\rho\colon\F_3^2\longrightarrow\Aut(E),
		\qquad \rho(a,b)(z):=\zeta_3^a z+bt.
	\]
	Put $m=n+2$ and identify $G\cong\F_3^n$.
	Let $\theta_i=(\lambda_i,\lambda_{i+1})$ be the maps constructed in the proof of Proposition~\ref{prop:p3-upper} for $e=1$.
	Define $\rho_i:=\rho\circ\theta_i\colon G\to\Aut(E)$ and let $G$ act diagonally on $E^m$ by
	\[
		h\cdot(z_i)_{i=1}^m:=\bigl(\rho_i(h)(z_i)\bigr)_{i=1}^m.
	\]
	Since $\rho$ is faithful, $\ker\rho_i=\ker\theta_i$, so~\hyperref[eq:p3-minimality]{\textup{(M$_3$)}} proves minimality.
    
	The fixed-point-free elements of $\rho(\F_3^2)$ are precisely the non-trivial translations $\rho(0,b)$.
	Hence
	\[
		\Sigma_\rho=\langle(1,0)\rangle
		\cup\langle(1,1)\rangle
		\cup\langle(1,-1)\rangle.
	\]
	The stabilizer set of $\rho_i$ is $\theta_i^{-1}(\Sigma_\rho)$, so~\hyperref[eq:p3-freeness]{\textup{(F$_3$)}} proves freeness.

	By \cite{BG21}*{Proposition~3.1}, infinitesimal rigidity of the diagonal action is equivalent to the absence of $G$-invariant quadratic differentials $dz_i\otimes dz_j$.
	Since
	\[
		\rho_i(h)^*dz_i=\zeta_3^{\lambda_i(h)}dz_i,
	\]
	it suffices to show that
	\[
		\lambda_i+\lambda_j\neq0
		\qquad\text{for all }1\le i,j\le m.
	\]
	The forms $\lambda_1,\ldots,\lambda_4$ are pairwise linearly independent, while each remaining form $\lambda_{4+s}$, $1\le s\le n-2$, contains the variable $W_s$, which occurs in no other form.
	Thus every $\lambda_i$ is non-zero and no two are negatives of one another, proving infinitesimal rigidity.
	Since the action is free, the quotient map $E^m\to E^m/G$ is unramified and
	\[
		H^1(E^m/G,\Theta_{E^m/G})
		\cong H^1(E^m,\Theta_{E^m})^G
		=0.
	\]
	Finally, $E^m$ has Kodaira dimension~$0$, and since the Kodaira dimension is invariant under unramified covers, the same holds true for $E^m/G$.
\end{proof}

\subsection{Abelian mixed actions} \label{subsec:mixed}

Having dropped the hyperbolicity condition in the preceding subsection, we finally analyse what happens if we consider mixed actions.
For this, we first recall some necessary results.
Every automorphism of a product of curves of genus at least two lifts to its universal cover $\Delta^n$, and
\[
	\Aut(\Delta^n)=\Aut(\Delta)^n\rtimes\mathfrak S_n.
\]
It therefore acts by factorwise automorphisms followed by a permutation of isomorphic factors.
By \cite{FGR26}*{Lemma~2.4}, a finite action on $C_1\times\cdots\times C_n$ induces a homomorphism $\pi\colon G\to\mathfrak S_n$.
Following \cite{FGR26}*{Definition~2.5}, the action is called \emph{unmixed} if $\pi$ is trivial and \emph{mixed} otherwise.
Thus mixed actions arise precisely from a non-trivial permutation part, and $G^0:=\ker(\pi)\trianglelefteq G$ is the subgroup preserving every factor.

We first recall why no mixed Beauville surface can have an abelian acting group.
In the surface case, a mixed action exchanges the two factors, which must therefore be isomorphic.
After identifying both with a curve $C$, the subgroup $G^0$ has index two, and one obtains an exact sequence
\[
	1\longrightarrow G^0\longrightarrow G \longrightarrow\Z/2\Z\longrightarrow1.
\]
In particular, $G$ contains an involution $\sigma$.

Suppose now that $G$ is abelian.
By the classification of abelian unmixed Beauville surface groups, there exists $n$ such that $G^0\cong(\Z/n\Z)^2$ with $\gcd(n,6)=1$, and thus necessarily $\sigma \not \in G^0$.
Thus there are automorphisms $\alpha,\beta\in\Aut(C)$ such that
\[
	\sigma(x,y)=(\alpha(y),\beta(x)).
\]
The relation $2 \sigma=0$ implies $\beta=\alpha^{-1}$, and hence
\[
	\operatorname{Fix}(\sigma)=\{(x,\beta(x))\mid x\in C\},
\]
contradicting freeness.

If we have more than two factors, then the preceding argument does not apply, and we can construct mixed actions of abelian groups:

Let $p\ge5$ be a prime and put
\[
	H=(\Z/p\Z)^2=\langle e_1,e_2\rangle.
\]
For example, the two triples
\[
	\begin{aligned}
		S_1&=(e_1,e_2,-e_1-e_2),\\
		S_2&=(e_1+2e_2,e_1-e_2,-2e_1-e_2)
	\end{aligned}
\]
form a classical Beauville structure on $H$.
Let $C_i\to\PP^1$ be the triangle $H$-cover associated with $S_i$.
Each $C_i$ has genus
\[
	g(C_i)=1+\frac{p(p-3)}2\geq2,
\]
and $H$ acts freely on $Z:=C_1\times C_2$.
Set
\[
	G=\langle\tau\rangle\times\langle e_2\rangle
	\cong \Z/p^2\Z \times \Z/p\Z,
	\qquad e_1=p\tau,
\]
so that $H=\langle e_1,e_2\rangle$ is an index-$p$ subgroup of $G$, and put
\[
	Y:=Z^p\cong C_1^p\times C_2^p.
\]
Writing points of $Y$ as $(z_0,\ldots,z_{p-1})$ with $z_j\in Z$, define
\[
	\begin{aligned}
		h\cdot(z_0,\ldots,z_{p-1})
		&=(h\cdot z_0,\ldots,h\cdot z_{p-1}) &&(\text{for } h\in H),\\
		\tau\cdot(z_0,\ldots,z_{p-1})
		&=(z_1,\ldots,z_{p-1},e_1\cdot z_0).
	\end{aligned}
\]
Then the $p$-fold iterate of $\tau$ acts diagonally as $e_1=p\tau$, and $\tau$ commutes with the diagonal $H$-action.
Hence these formulas define a faithful action of $G$ on $Y$: an element outside $H$ has non-trivial permutation part, whereas the restriction to $H$ is faithful.
The action is not diagonal, since $\tau$ cyclically permutes the $p$ copies of $Z$.

The action is free: on the one hand, $H$ acts freely on $Z$, so its diagonal action on $Y$ is free.
On the other hand, if $g\in G\setminus H$, then $g=k\tau+h$ for some $h\in H$ and $1\leq k\leq p-1$.
Since $G$ is abelian and $H$ has exponent $p$,
\[
	pg=kp\tau+ph=ke_1\neq0.
\]
If $g$ fixed a point of $Y$, then $pg=ke_1$ would fix that point as well, contradicting the freeness of the $H$-action.

Moreover, $X_p := Y/G$ is rigid.
Since the $C_i$ are triangle curves,
\[
	H^1(C_i,\Theta_{C_i})^H=0.
\]
The K\"unneth decomposition therefore gives
\[
	H^1(Y,\Theta_Y)^G
	\subseteq H^1(Y,\Theta_Y)^H=0,
\]
and hence $H^1(X_p,\Theta_{X_p})=0$.
Consequently, $X_p$ is a rigid quotient of a product of $2p$ curves by a free mixed action of the abelian group $G$.
Notice also that Theorem~\ref{thm:main} gives $\du(G)=\infty$, because the largest invariant factor $p^2$ has multiplicity one.
Thus $G$ admits no unmixed Beauville structure.

We obtain the following corollary by the same construction, using Proposition~\ref{prop:EllipticRigid} in the elliptic case.

\begin{corollary}\label{cor:abelian-mixed-realizations}
	Let $G$ be a non-trivial finite abelian group of odd order with non-cyclic Sylow subgroups.
	Then $G$ admits a free action on a product of curves of genus at least one with rigid quotient.
	Such an action can be chosen unmixed if and only if $m_{\exp}(G)\ge2$.
\end{corollary}
\begin{proof}
	Write $G=\langle\tau\rangle\oplus B$, where $\ord(\tau)=\exp(G)=n_1$ and $\exp(B)=n_2\mid n_1$, and put
	\[
		k:=\frac{n_1}{n_2},
		\qquad H:=\langle k\tau\rangle\oplus B.
	\]
	Then $G/H\cong\Z/k\Z$ and $m_{\exp}(H)\ge2$.
	Moreover, $G[p]\subseteq H$ for every $p\mid|G|$, since the Sylow subgroups are non-cyclic.
	Choose the action $H\curvearrowright Z$ from Theorem~\ref{thm:main} if $\exp(H)>3$, and from Proposition~\ref{prop:EllipticRigid} if $\exp(H)=3$.
	As above, extend the diagonal $H$-action on $Y:=Z^k$ to $G$ by
	\[
		\tau\cdot(z_0,\ldots,z_{k-1})
		=(z_1,\ldots,z_{k-1},(k\tau)\cdot z_0).
	\]
	Every $0\ne g\in G$ has a non-zero multiple in some $G[p]\subseteq H$, so the action is free.

	If $\exp(H)>3$, K\"unneth's formula gives $H^1(Y,\Theta_Y)^H=0$.
	In the elliptic case, repeating the factors preserves $\lambda_i+\lambda_j\ne0$, so the same vanishing follows from \cite{BG21}*{Proposition~3.1}.
	Hence $H^1(Y,\Theta_Y)^G=0$, and the quotient is rigid.

	If $m_{\exp}(G)\ge2$, then $H=G$ and the action is unmixed.
\end{proof}

\begin{remark}
	The preceding examples show that a group $G$ with $\du(G)=\infty$ can still admit a free action on a product of curves with rigid quotient if elliptic factors or mixed actions are allowed.
\end{remark}

\begingroup
	\renewcommand{\addcontentsline}[3]{}
	\section*{Declaration of generative AI assistance}

	The author used \href{https://openai.com/codex/}{OpenAI Codex} with the GPT-5.6 Sol model (August 2026) to help improve the formulation of selected passages.
	The author reviewed all AI-assisted content and takes full responsibility for the final manuscript.
\endgroup


\begin{thebibliography}{xxxxxxxx}

	\bib{BBF12}{article}{
		author={Barker, Nathan W.},
		author={Boston, Nigel},
		author={Fairbairn, Ben T.},
		title={A note on Beauville $p$-groups},
		journal={Exp. Math.},
		volume={21},
		date={2012},
		number={3},
		pages={298--306},
		doi={10.1080/10586458.2012.669267}
	}

	\bib{BCG05}{incollection}{
		author={Bauer, Ingrid},
		author={Catanese, Fabrizio},
		author={Grunewald, Fritz},
		title={Beauville surfaces without real structures},
		booktitle={Geometric methods in algebra and number theory},
		editor={Bogomolov, Fedor},
		editor={Tschinkel, Yuri},
		series={Progr. Math.},
		volume={235},
		publisher={Birkh\"auser Boston},
		address={Boston, MA},
		date={2005},
		pages={1--42},
		doi={10.1007/0-8176-4417-2\_1}
	}

	\bib{BG21}{article}{
		author={Bauer, Ingrid},
		author={Gleissner, Christian},
		title={Towards a classification of rigid product quotient varieties of Kodaira dimension $0$},
		journal={Boll. Unione Mat. Ital.},
		volume={15},
		date={2022},
		number={1-2},
		pages={17--41},
		doi={10.1007/s40574-021-00295-4}
	}

    \bib{B83}{book}{
		author={Beauville, Arnaud},
		title={Complex algebraic surfaces},
		series={London Math. Soc. Lecture Note Ser.},
		volume={68},
		note={Translated from the French by R. Barlow, N. I. Shepherd-Barron and M. Reid},
		publisher={Cambridge Univ. Press},
		address={Cambridge},
		date={1983}
	}

	\bib{Cat00}{article}{
		author={Catanese, Fabrizio},
		title={Fibred surfaces, varieties isogenous to a product and related moduli spaces},
		journal={Amer. J. Math.},
		volume={122},
		date={2000},
		number={1},
		pages={1--44},
		doi={10.1353/ajm.2000.0002}
	}

	\bib{CF22}{article}{
		author={Carta, Ludo},
		author={Fairbairn, Ben T.},
		title={Generalised Beauville groups},
		journal={J. Group Theory},
		volume={25},
		date={2022},
		number={2},
		pages={355--377},
		doi={10.1515/jgth-2021-0034}
	}

	\bib{FAG17}{article}{
		author={Fern\'andez-Alcober, Gustavo A.},
		author={G\"ul, \c{S}\"ukran},
		title={Beauville structures in finite $p$-groups},
		journal={J. Algebra},
		volume={474},
		date={2017},
		pages={1--23},
		doi={10.1016/j.jalgebra.2016.11.007}
	}

	\bib{FGGS24}{article}{
		author={Fern\'andez-Alcober, Gustavo A.},
		author={Gavioli, Norberto},
		author={G\"ul, \c{S}\"ukran},
		author={Scoppola, Carlo M.},
		title={Beauville $p$-groups of wild type and groups of maximal class},
		journal={Monatsh. Math.},
		volume={205},
		date={2024},
		number={4},
		pages={725--734},
		doi={10.1007/s00605-024-01982-y}
	}

	\bib{Fai19}{incollection}{
		author={Fairbairn, Ben T.},
		title={Beauville $p$-groups: a survey},
		booktitle={Groups St Andrews 2017 in Birmingham},
		series={London Math. Soc. Lecture Note Ser.},
		volume={455},
		publisher={Cambridge Univ. Press},
		date={2019},
		pages={271--288},
		doi={10.1017/9781108692397.012}
	}

	\bib{FG16}{article}{
		author={Frapporti, Davide},
		author={Gleissner, Christian},
		title={On threefolds isogenous to a product of curves},
		journal={J. Algebra},
		volume={465},
		date={2016},
		pages={170--189},
		doi={10.1016/j.jalgebra.2016.06.034}
	}

	\bib{FGR26}{article}{
		author={Fallucca, Federico},
		author={Gleissner, Christian},
		author={Ruhland, Noah},
		title={On rigid varieties isogenous to a product of curves},
		journal={J. Algebra},
		volume={688},
		date={2026},
		pages={393--419},
		doi={10.1016/j.jalgebra.2025.09.016}
	}

	\bib{FJ11}{article}{
		author={Fuertes, Yolanda},
		author={Jones, Gareth A.},
		title={Beauville surfaces and finite groups},
		journal={J. Algebra},
		volume={340},
		date={2011},
		number={1},
		pages={13--27},
		doi={10.1016/j.jalgebra.2011.04.027}
	}

	\bib{FMP13}{article}{
		author={Fairbairn, Ben T.},
		author={Magaard, Kay},
		author={Parker, Christopher},
		title={Generation of finite quasisimple groups with an application to groups acting on Beauville surfaces},
		journal={Proc. Lond. Math. Soc.},
		volume={107},
		date={2013},
		number={4},
		pages={744--798},
		doi={10.1112/plms/pds097}
	}

	\bib{GT12}{incollection}{
		author={Gonz\'alez-Diez, Gabino},
		author={Torres-Teigell, David},
		title={An introduction to Beauville surfaces via uniformization},
		booktitle={Quasiconformal mappings, Riemann surfaces, and Teichm\"uller spaces},
		series={Contemp. Math.},
		volume={575},
		publisher={American Mathematical Society},
		address={Providence, RI},
        doi={10.1090/conm/575/11392},
		date={2012},
		pages={123--151}
	}

	\bib{GLL12}{article}{
		author={Garion, Shelly},
		author={Larsen, Michael},
		author={Lubotzky, Alexander},
		title={Beauville surfaces and finite simple groups},
		journal={J. Reine Angew. Math.},
		volume={666},
		date={2012},
		pages={225--243},
		doi={10.1515/CRELLE.2011.117}
	}

	\bib{GM12}{article}{
		author={Guralnick, Robert},
		author={Malle, Gunter},
		title={Simple groups admit Beauville structures},
		journal={J. Lond. Math. Soc.},
		volume={85},
		date={2012},
		number={3},
		pages={694--721},
		doi={10.1112/jlms/jdr062}
	}

	\bib{Jam77}{article}{
		author={Jamison, Robert E.},
		title={Covering finite fields with cosets of subspaces},
		journal={J. Combin. Theory Ser. A},
        doi={10.1016/0097-3165(77)90001-2},
		volume={22},
		date={1977},
		number={3},
		pages={253--266}
	}

	\bib{LM02}{book}{
		author={Leedham-Green, Charles R.},
		author={McKay, Susan},
		title={The structure of groups of prime power order},
		series={London Mathematical Society Monographs. New Series},
		volume={27},
		publisher={Oxford University Press},
		address={Oxford},
		date={2002},
		doi={10.1093/oso/9780198535485.001.0001}
	}

\end{thebibliography}
\end{document}